\documentclass[journal, twocolumn, 10pt, letterpaper]{IEEEtran}
\usepackage{amsmath,amssymb,array,graphicx,verbatim}

\usepackage{latexsym,url}
\usepackage{amsmath,amssymb,amsthm}
\usepackage{algorithm,algorithmic}
\usepackage{color,array,graphicx,verbatim,xtab}
\usepackage{cite} % introduces sort and compress of multiple refs
\usepackage{setspace}

\usepackage
  [breaklinks,bookmarks,bookmarksnumbered,bookmarksopen,bookmarksopenlevel=2]
  {hyperref}

\newtheorem{remark}{Remark}

\usepackage{color}
\newcommand{\mycomment}[3]%
{%
\marginpar{%
  \hfil%
  \tiny{\textcolor{#2}{{\bf\textsc{#1}}}}%
  \hfil%
}%
\footnote{\textcolor{#2}{{\bf\textsc{#1}:}~~#3}}%}
}
\newcommand{\kakhbod}[1]%
{\mycomment{kakhbod}{blue}{#1}}
\newcommand{\jckoo}[1]%
{\mycomment{jckoo}{red}{#1}}

\begin{document}

%\title{Market Mechanisms with Non-Price-Taking Agents}
  \title{Correction to ``An Efficient Game Form for Unicast Service Provisioning"}
%\author
%{\Large{Ali Kakhbod
% and Demosthenis Teneketzis~\IEEEmembership{\Large Fellow, IEEE}
% }%
% \thanks{Ali Kakhbod and Demosthenis Teneketzis are with the department of
% EECS, University of Michigan, Ann Arbor, MI -- 48109, USA. (email:
% \texttt{\{akakhbod,teneket\}@umich.edu}). Their work was supported in part
% by NSF Grant CCR-0325571 and by NSF Grant CCF-1111061.}}
%\textsc{preliminary paper (please, do not circulate or post online)}} 

\author{
{\Large{Ali Kakhbod}}\thanks{A. Kakhbod is with the Department of Electrical and System Engineering University of Pennsylvania, Philadelphia, (email:
\texttt{akakhbod@seas.upenn.edu}).} \quad  {\Large{and}} \quad 
{\Large{Demosthenis Teneketzis}} \thanks{D. Teneketzis is with the Department of Electrical Engineering and Computer Science, University of Michigan (email:
\texttt{teneket@eecs.umich.edu}).}
\\
{May 2013}}

\maketitle

\begin{abstract}
%This note corrects an error in \textbf{\cite{kakh}}.
A correction to the specification of the mechanism proposed in \textbf{\cite{kakh}} is given.

\end{abstract}

%Due to an error,  t
The mechanism presented in \textbf{\cite{kakh}} has a tax function which is not differentiable with respect to the allocations. We need a tax function which is differentiable with respect to the allocations so that we can have Nash implementation. We correct it as follows. 

We consider the problem formulated in \textbf{\cite{kakh}}. We use the same notation as in \textbf{\cite{kakh}}.  

%We describe the new mechanism and analyze its properties.  

\noindent{\textbf{Specification of the game form/mechanism}:
%\label{main}
%\section{Specification of the game form/mechanism}\label{main}

\noindent{\emph{Message space}}: The message space is the same as that of the mechanism presented in \textbf{\cite{kakh}}. 
A message of  user $i\in \mathcal{N}$ ($\mathcal{N}$ denotes the set of users) is of the form 
\begin{align*}
m_i=({x}_i,p_i^{l_{i_1}},p_i^{l_{i_2}},\cdots,p_i^{l_{i_{|\mathcal{R}_i|}}}),
\end{align*}
where ${x}_i$ denotes the (non-negative) bandwidth user $i$ requests at all the links of his route,  and $p_i^{l_{i_k}}\geq 0$ denotes the price user $i$ is willing to pay per unit of bandwidth at link $l_{j_k}$ of his route $\mathcal{R}_i$.
%\medskip

\noindent{\emph{Outcome function}}:
For any $m\in \mathcal{M}$, the outcome function is defined as follows:
\begin{align*}
f(m)&=(x_1,x_2,\cdots,x_n,t_1,t_2,\cdots,t_n) \\
%\end{align*}
%%where, for any $i\in \mathcal{N}$:
% %\begin{align}\label{test}
% %x_i= \hat{x}_i  \min\{\frac{c^l}{\sum_{j\in \mathcal{G}^l}\hat{x}_j}: l\in \mathcal{R}_i\} \ \  \mathbb{1}\{\hat{x}_i>0\},
% %\end{align}
%and 
% \begin{align*}
 t_i&=\sum_{l\in \mathcal{R}_i} t_i^l,
 \end{align*}
where $t_i^l$ is the tax paid by user $i$ for using link $l$. The form of $t_i^l$ is the same as the tax function defined in \textbf{\cite{kakh}} excluding the term that is of the form described by relation  \textbf{(23) in \cite{kakh}}. For example, if $|\mathcal{G}^l|>3$,  ($\mathcal{G}^l$ denotes the set of users using link $l$) the tax function in Eq. \textbf{(13) of \cite{kakh}}  now becomes,
% 
% 
%where $c^l$ is the capacity of link $l$, $\mathcal{G}^l$ denotes the set of users using link $l$, and $t_i^l$ is  the same as the tax function defined in \textbf{\cite{kakh}} excluding the term that is of the form described by relation  \textbf{(23) in \cite{kakh}}.
%For example, if $|\mathcal{G}^l|>3$, the tax function in Eq. \textbf{(13) of \cite{kakh}}  now becomes
\begin{align}\label{tax}
t_i^l&=P_{-i}^lx_i+(p_i^l-P_{-i}^l-\zeta_+^l)^2 \nonumber \\
&\quad -2P_{-i}^l\left(p_i^l-P_{-i}^l \right)\left(\frac{\mathcal{E}_{-i}^l+x_i}{\gamma}\right)+\Phi_{i}^l,
\end{align} 
where
\begin{align}
\zeta_+^l=\max\{0, \frac{\sum_{i\in \mathcal{G}^l}x_i-c^l}{\hat{\gamma}}\},
\end{align}
$c^l$ is the capacity of link $l$,  $\Phi_i^l$ is defined by Eq. \textbf{(14) in \cite{kakh}}, 
\begin{align}
P_{-i}^l=\frac{\sum_{\substack{j\in \mathcal{G}^l\\ j\neq i}}p_j^l}{|\mathcal{G}^l|-1},   \quad \quad \quad
 \mathcal{E}_{-i}^l=\sum_{\substack{j\in \mathcal{G}^l\\ j\neq i}}x_j-c^l,
\end{align}
($P_{-i}^l$ and $ \mathcal{E}_{-i}^l$ are the same as in \textbf{\cite{kakh}}) and $\gamma, \hat{\gamma}$,  are positive constants.

This completes the specification of the mechanism. 

Based on the above specification, the proof of Lemma \textbf{2 in \cite{kakh}} is updated as follows.

%\section{Properties of a NE}
\begin{proof}[\textbf{Proof of Lemma 2 in \cite{kakh}}]
%\noindent{\textbf{Proof of Lemma 2 in \cite{kakh}}}:
Let $m^*=(m_i^*,m_{-i}^*)$ be a NE of the game induced by the mechanism. Since user $i$ does not control $\Phi_i^l$, it implies 
$\frac{\partial \Phi_i^l}{\partial p_i^l}=0$, (as in Eq. \textbf{(34)} \textbf{of \cite{kakh}}).
By following the same steps as in equations \textbf{(35-38) of \cite{kakh}}, we obtain for any $l\in \textbf{L}$:
\begin{align}\label{s1}
\frac{\partial t_i^l}{\partial p_i^l}\big|_{m=m^*}\!\!=2\left[(p_i^{*l}-P_{-i}^{*l}-\zeta_+^{*l})-P_{-i}\left(\frac{\mathcal{E}_{-i}^{*l}+x_i^*}{\gamma}\right)\right]=0.
\end{align}
Summing \eqref{s1} over all $i\in \mathcal{G}^l$, we get
\begin{align}
\sum_{i\in \mathcal{G}^l}\frac{\partial t_i^l}{\partial p_i^l}\big|_{m=m^*}&=\!\!\sum_{i\in \mathcal{G}^l}\left[(p_i^{*l}-P_{-i}^{*l}-\zeta_+^{*l})-P_{-i}\left(\frac{\mathcal{E}_{-i}^{*l}+x_i^*}{\gamma}\right)\right] \nonumber\\ 
&=-|\mathcal{G}^l|\zeta_+^{*l}-\sum_{i\in \mathcal{G}^l}P_{-i}^{*l}\left(\frac{\mathcal{E}_{-i}^{*l}+x_i^*}{\gamma}\right)   \nonumber\\ \label{s2}
&=0.
\end{align}
Suppose $\sum_{i\in \mathcal{G}^l}x_i^*>c^l$. Then we must have, $\zeta_+^{*l}>0$ and $\sum_{i\in \mathcal{G}^l}P_{-i}^{*l}\left(\frac{\mathcal{E}_{-i}^{*l}+x_i^*}{\gamma}\right)\geq 0$. But this contradicts Eq. \eqref{s2}.  Therefore, we must have
\begin{align}\label{xx}
\sum_{i\in \mathcal{G}^l}x_i^*\leq c^l.
\end{align}
This implies, 
%since $\gamma$ is strictly positive and $P_{-i}^{*l}\geq 0$, we conclude from \eqref{s2} that we must have  
\begin{align}\label{ss}
\zeta_+^{*l}=0.
\end{align}
Combining \eqref{ss} along with \eqref{s2} we obtain
\begin{align}\label{tz}
\sum_{i\in \mathcal{G}^l}P_{-i}^{*l}\left(\frac{\mathcal{E}_{-i}^{*l}+x_i^*}{\gamma}\right)=0.
\end{align}
Moreover, combining \eqref{xx} and \eqref{tz} we obtain
\begin{align}\label{s4}
P_{-i}^{*l}\left(\frac{\mathcal{E}_{-i}^{*l}+x_i^*}{\gamma}\right)=0.
\end{align}
for every $i\in \mathcal{G}^l$.
Using  \eqref{ss} and  \eqref{s4} in \eqref{s1} we obtain
\begin{align}\label{xc}
p_i^{*l}=P_{-i}^{*l}.
\end{align}
Since \eqref{xc} is true for all $i\in \mathcal{G}^l$, it implies,
\begin{align}
p_i^{*l}=p_j^{*l}=P_{-i}^{*l}&=:p^{*l}, \label{t4}
\end{align}
and along with \eqref{s4} it implies
\begin{align}
p^{*l} \mathcal{E}^{*l}=0, \label{t5}
\end{align}
where $\mathcal{E}^{*l}=\sum_{i\in \mathcal{G}^l}x_i^*-c^l$  ($\mathcal{E}^{*l}$ is the same as in \textbf{\cite{kakh}}).

Furthermore, since
\begin{align}\label{cc}
 \frac{\partial \Phi_i^l}{\partial x_i}=0  
\end{align}  
 (Eq. \textbf{(34)} in \textbf{\cite{kakh}})), it follows from \eqref{tax} that
\begin{align}\label{x1}
\frac{\partial t_i^l}{\partial x_i}\big|_{m=m^*}=p^{*l}.
\end{align}
because  of \eqref{ss}, \eqref{t4}, \eqref{t5}, and \eqref{cc}. 
%Note that if $x_i^*=0$, since user $i$ does not have incentive to increase its demand, it follows that $\frac{\partial t_i^l}{\partial x_i}\big|_{m=m^*}\leq p^{*l}$.
\end{proof}
\begin{remark}
\emph{The proof of Theorem \textbf{5} follows when $x_i^*>0$. Note that, when $x_i^*=0$, since user $i$ does not have incentive to increase its demand, it follows that 
\begin{align}\label{de}
\frac{\partial \textbf{U}_i(x_i)}{\partial x_i}-\sum_{l\in \mathcal{R}_i} p^{*l}\big|_{m=m^*}\leq 0.
\end{align}
Now,  set $\lambda^{*l}=p^{*l}$. Then \eqref{t5} and \eqref{de} are consistent with the KKT conditions \textbf{(68-70) of \cite{kakh}}.  }
\end{remark}

\section{Properties of the mechanism}
\noindent{\textbf{Existence of Nash equilibria (NE)}}:  The proof of existence of NE of the game induced by the mechanism  is the same as in\textbf{\cite{kakh}} (see Theorem \textbf{6}, page  \textbf{398}, and its proof in \textbf{\cite{kakh}}; also see the proof of Theorem \textbf{7}).

\noindent{\textbf{Feasibility of allocations at NE}}: Because  of the specification of the mechanism and  Eq. \eqref{ss},  the allocations corresponding to all NE are in the feasible set.

\noindent{\textbf{Budget Balance at any feasible allocation}}: Budget balance at any feasible allocation  follows by Lemma \textbf{3} of \textbf{\cite{kakh}}.

\noindent{\textbf{Individual Rationality}}: Individual rationality follows by Theorem \textbf{4} of \textbf{\cite{kakh}}.

\noindent{\textbf{Nash implementation}}: 
Nash implementation follows by Theorem \textbf{5} of \textbf{\cite{kakh}}.

\medskip
\noindent{\emph{Acknowledgment}}: The authors wish to thank Professor A. Anastasopoulos for  pointing out the error.
%pointing out to them the error in equation (64) of the paper. 

\bibliographystyle{nonumber}

\end{document}